\documentclass{article}
\usepackage[utf8]{inputenc}
\usepackage{float}
  \usepackage{amsmath}
  \usepackage{amsthm}
 \usepackage{amssymb}
\usepackage{amsfonts}
\usepackage{graphics}
\usepackage{bbm}
 \usepackage[english]{babel}
\usepackage{hyperref}
 \newlength\tindent
\newtheorem{theorem}{Theorem}

\newtheorem{lemma}{Lemma}

\newtheorem{definition}{Definition}

\newtheorem{remark}{Remark}
\newtheorem{problem}{Problem}

\begin{document}
\title{Klarner's sequences have density $0$}
\author{Alexander Kalmynin \footnote{National Research University Higher School of Economics, Moscow, Russia} \\ email: \href{mailto:alkalb1995cd@mail.ru}{alkalb1995cd@mail.ru}}
\date{}
\maketitle
\begin{abstract}
We prove that the set of non-negative integers $n$ which can be obtained from $0$ by repeated application of affine maps $x\mapsto 2x$, $x\mapsto 3x+2$ and $x\mapsto 6x+3$ has density $0$. Moreover, we reach the same conclusion for five other sequences of similar type introduced by D. Klarner. More generally, we establish
a criterion excluding stationary measures with bounded density on the profinite group $\mathbb Z_P\cong\prod_{p\mid P}\mathbb Z_p$.
A density-to-measure construction using invariant Banach limits and
a combinatorial argument involving rational orbits and affine words
yield our density-zero results. The criterion also applies to an infinite
family of mixed-slope affine semigroups.
\end{abstract}

\section{Introduction and main results}

Let $A_1,\ldots,A_r$ be affine maps $A_i:x\mapsto a_ix+b_i$. We will assume that $a_i\in \mathbb Z_{>0}$ and $b_i\in \mathbb Z_{\geq 0}$, so that each $A_i$ maps non-negative integers to non-negative integers. Let $W=W(A_1,\ldots,A_r)$ be the semigroup generated by these maps (for convenience, we assume that $id\in W$). We study its orbits on sets of non-negative integers. More precisely, for $B\subseteq \mathbb Z_{\geq 0}$, we write $W(B)$ for the set of all integers of the form $w(b)$ with $w\in W$ and $b\in B$. The arithmetic structure of such orbit sets has been investigated since the work
of Klarner and Rado \cite{KlarnerRado1974}. A basic question, already raised in this circle of
problems by Erd\H{o}s and Graham \cite{ErdosGraham1980}, is the following

\begin{problem}
Given $W$ and $B\subset \mathbb Z_{\geq 0}$, under what conditions does $W(B)$ have zero asymptotic density?
\end{problem}

One particularly persistent instance arises from the three maps
\[
 x\longmapsto 2x,\qquad x\longmapsto 3x+2,\qquad
 x\longmapsto 6x+3.
\]
The question whether the orbit of $0$ has positive density is among the
problems of Klarner publicized by Guy \cite{Guy1983,Guy2004}; see also
Lagarias's historical and mathematical account \cite{Lagarias2016}.
Klarner identified six triples of translations for the slopes $2,3,6$
that generate free semigroups \cite{Klarner1982}. The same six examples
have subsequently featured in work on the geometry and growth of affine
semigroups, notably in Kolpakov--Talambutsa \cite{KolpakovTalambutsa2022}
and Shamazov--Talambutsa \cite{ShamazovTalambutsa2026}.
We prove that all six orbit sets starting at $0$ have density zero.

The role of the slopes is essential. For a finite family with
$\sum\limits_{i=1}^r a_i^{-1}<1$, classical counting arguments of Erd\H{o}s,
Klarner and Rado give sublinear bounds for the orbit of a finite set
\cite{KlarnerRado1974,Lagarias2016}. The borderline case
$\sum_{i=1}^r a_i^{-1}=1$ is substantially more delicate.
Freeness of the semigroup does not by itself resolve the density question:
the results of Shamazov and Talambutsa yield, in the balanced free case,
the almost-linear lower bound
\[
 \#\left(W(\{b\})\cap[0,x]\right)
 \gg \frac{x}{(\log x)^{(r-1)/2}}
\]
for a fixed starting point $b$, but this bound does not imply positive
density \cite{ShamazovTalambutsa2026}. Our result provides explicit
balanced free semigroups whose orbits nevertheless have density zero.
For comparison, the equal-slope case has been studied by Klarner
\cite{Klarner1981} and, more recently, by Konyagin, Protasov and
Talambutsa using refinement equations and Fourier analysis
\cite{KonyaginProtasovTalambutsa2025}. The mixed-slope setting considered
here calls for a different approach.

\begin{definition}
A set of generators $A_1,\ldots,A_r$ of a semigroup $W$ is called balanced if
\[
\sum_{i=1}^r\frac{1}{a_i}=1.
\]
\end{definition}

Next, let $P$ be the product of all primes $p$ which divide $\prod\limits_{i\leq r}a_i$. Consider the space $\mathbb Z_P$ of all $P$-adic integers: it can be viewed as an inverse limit of $\mathbb Z/P^N\mathbb Z$ with respect to projection maps $\mathbb Z/P^{N+1}\mathbb Z\to \mathbb Z/P^N\mathbb Z$ given by $x\mapsto x \pmod{P^N}$. This is a compact group and we also have an isomorphism
\[
\mathbb Z_P\cong \prod_{p\mid P}\mathbb Z_p.
\]
The Haar measure $\mu$ on $\mathbb Z_P$ is defined by
\[
\mu(a+P^N\mathbb Z_P)=P^{-N}.
\]
For any measurable set $C\subset \mathbb Z_P$ and any affine map $\varphi(x)=cx+d$ with $c\in \mathbb Z_P\cap \mathbb Q_P^{*}$ we have $\mu(\varphi(C))=|c|_P\mu(C)$, where $|c|_P=\prod\limits_{p\mid P}|c|_p$ is the product of $p$-adic norms of $c$. In particular, when $\varphi \in W$ we always have $|c|_P=\frac{1}{c}$.

\begin{definition}
For a semigroup $W$ with a fixed set of generators $A_1,\ldots,A_r$ we call a measure $\nu$ on $\mathbb Z_P$ stationary for $W$ if
\[
\nu=\sum_{i=1}^r \frac{1}{a_i}(A_i)_*\nu.
\]
Here $(A_i)_*\nu(C)=\nu(A_i^{-1}C)$.
\end{definition}
Note that the definition of $W$-stationarity depends not only on the semigroup but also on the choice of generators. It is easy to see that if the generators are imbalanced, there are no finite non-zero $W$-stationary measures. Indeed, we would have
\[
\nu(\mathbb Z_P)=\sum_{i=1}^r \frac{1}{a_i}(A_i)_*\nu(\mathbb Z_P)=\nu(\mathbb Z_P)\sum_{i=1}^r \frac{1}{a_i},
\]
which implies $\nu(\mathbb Z_P)=0$ and $\nu=0$.

It turns out that a negative answer to Problem~1 in the balanced case
is intimately related to the existence of certain $W$-stationary
measures on $\mathbb{Z}_P$.

\begin{theorem}
Let $W=\langle A_1,\ldots,A_r\rangle$ be a semigroup with a
balanced set of generators. If there is a set $B\subseteq\mathbb{N}$
of density zero such that the lower density of $W(B)$ is nonzero,
then there is a measure $\nu\neq 0$ on $\mathbb{Z}_P$ which is
$W$-stationary and majorized by the Haar measure $\mu$:
\[
    \nu(A)\leq\mu(A)
    \qquad\text{for all measurable sets }A.
\]
\end{theorem}

We prove Theorem~1 in Section~2. The proof relies on the existence of
certain Banach limits on $\mathbb{N}$ with special invariance
properties. Note that for a general set $B$ Theorem~1 tells us nothing about the upper density.

Thus, to establish that $W(B)$ has zero lower density for every zero-density set $B$, it suffices to rule out the existence of such stationary measures.

\begin{definition}
We call a measure $\nu$ on $\mathbb{Z}_P$ a
\emph{KR-measure} (for $W$) if it is $W$-stationary and has the form
\[
    \nu=f\mu,
    \qquad f\in L^\infty(\mathbb{Z}_P,\mu).
\]
\end{definition}

By the Radon--Nikodym theorem, this is equivalent to
\[
    \nu(A)\leq \lVert f\rVert_\infty\,\mu(A)
\]
for all measurable sets $A$ (the Radon--Nikodym derivative
$d\nu/d\mu$ is essentially bounded).

 Our main criterion rules out nonzero KR-measures by combining a
real-dynamical covering condition with a finite collection of rational
orbits. The real dynamics is closely related to the geometric viewpoint
on freeness developed by Kolpakov and Talambutsa
\cite{KolpakovTalambutsa2022}; the arithmetic part of the argument is
instead carried out on $\mathbb Z_P$.

\begin{theorem}
Let $W=\langle A_1,\ldots,A_r\rangle$ be an affine semigroup
generated by a balanced set of generators $A_1,\ldots,A_r$.
Let $P$ be as above. Suppose that the following properties hold:
\begin{enumerate}
\item \emph{(Arithmetic condition.)}
$P=pq$ is a product of two distinct primes,
\[
    A_1x=p^N x
    \qquad\text{for some }N\in\mathbb{Z}_{>0},
\]
and the leading coefficient of $A_i x$ is divisible by $q$
for every $1<i\leq r$.

\item \emph{(Real dynamics condition.)}
There is a compact set $K\subset\mathbb{R}$ with nonempty interior
and a map $F:K\to K$ such that for every $x\in K$ there exists
$i\leq r$ with $F(x)=A_i x$.

\item \emph{(Discrete dynamics condition.)}
Let $v=\max_{i} v_q(a_{A_i})$ --- the maximal $q$-adic valuation among the leading coefficients of the generators. For every $s\in K\cap \mathbb Q$  of the form $s=a/q^j$ with $q\nmid a$ and $1\leq j\leq v$ there is an element $w\in W$ such that $ws=0$.
\end{enumerate}
Then there is no non-zero KR-measure for $W$. In particular, $W(B)$ has
lower density zero if $B$ has density zero.
\end{theorem}

\begin{remark}
Condition~(A) is chosen to model Klarner's semigroups with the leading
coefficients $2,3,6$.

Condition~(B) will usually be applied to the case when $K$ is a
closed interval and
\[
    K=\bigcup_{i\leq r}A_i^{-1}(K),
\]
where the preimages have disjoint interiors. We will see in further
sections that this is quite a restrictive condition. In fact, this
is Zariski closed in a suitable sense if we fix the leading
coefficients.

Condition~(C) involves a finite check, since the corresponding set of rational numbers is always finite: it is contained in the finite set $q^{-v}(q^vK\cap \mathbb Z)$.
\end{remark}

As an application of Theorem~2, we establish the following.

\begin{theorem}
Let $W$ be one of the six semigroups
\[
\begin{gathered}
\langle 2x,\,3x+2,\,6x+3\rangle,\qquad
\langle 2x,\,3x+3,\,6x+10\rangle,\\[3pt]
\langle 2x+2,\,3x,\,6x+15\rangle,\qquad
\langle 2x+1,\,3x,\,6x+2\rangle,\\[3pt]
\langle 2x+2,\,3x+1,\,6x\rangle,\qquad
\langle 2x+1,\,3x+6,\,6x\rangle.
\end{gathered}
\]
Then there is no non-zero KR-measure for $W$, and $W(B)$ has lower density zero
for every $B\subseteq\mathbb{N}$ of density zero.
\end{theorem}
To deduce that Klarner's original sequence has density zero, rather than merely zero lower density, we prove the following result.
\begin{theorem}
 Suppose that $W$ is balanced. Let $\overline{d}$ and $\underline{d}$ be the upper and lower densities of $W(B)$ respectively. Suppose that there is a constant $c>0$ such that the counting function $b(x)$ of $B$ satisfies $b(x)=O(x^{1-c})$. Then
\[
\overline{d}\leq \max{a_i}\cdot\underline{d}.
\]
In particular, if $\underline{d}=0$, then $W(B)$ has density $0$.
\end{theorem}

Section~3 will be devoted to the proof of Theorem~2. In Section~4
we check that Theorem~2 applies (up to some changes of variables) to the six semigroups from
Theorem~3.

Finally, Section~5 discusses further examples and the relation of our approach to Fourier-analytic methods.

\section{From sets to measures}

Here we prove Theorem~1. To do so, we use a nice class of Banach limits.

\begin{definition}
A functional $f:\ell^\infty\to \mathbb C$ is called a strongly invariant Banach limit if
\begin{itemize}
    \item[A)] $f1=1$, where $1=(1,1,\ldots)\in \ell^\infty$ --- the unit of $\ell^\infty$.
    \item[B)] $fx\geq 0$ for all $x\geq 0$.
    \item[C)] $fTx=fx$, where $Tx=(0,x_1,x_2,\ldots)$ --- the right shift operator.
    \item[D)] For all $m\in \mathbb N$ we have
    \[
    f\sigma_mx=fx,
    \]
    where $\sigma_mx=(\underbrace{x_1,\ldots,x_1,}_{m\text{ times}}\underbrace{x_2,\ldots,x_2}_{m\text{ times}},\ldots)$.
    \item[E)] $fCx=fx$, where $C$ is the Ces\`aro averaging operator
    \[
    Cx=\left(x_1,\frac{x_1+x_2}{2},\frac{x_1+x_2+x_3}{3},\ldots\right)
    \]
\end{itemize}
\end{definition}
\begin{remark}
Properties C) and D) together imply that
\[
f\pi_mx=\frac{1}{m}fx,
\]
where $\pi_mx=(x_1,\underbrace{0,\ldots,0,}_{(m-1)\text{ zeros}}x_2,\ldots)$.
\end{remark}
\begin{lemma}
Strongly invariant Banach limits exist.
\end{lemma}
\begin{proof}
See \cite{SemenovEtAl2020a, SemenovEtAl2020b}.
\end{proof}
From now on, let us fix a strongly invariant Banach limit $f$.
\begin{lemma}
Let $P$ be a squarefree integer and $A\subseteq \mathbb N$ be a subset. Define a measure $\nu=\nu_{A,f}$ on the algebra $\mathcal C$ of subsets of $\mathbb Z_P$ generated by $a+P^m\mathbb Z_P$ ($a\in \mathbb Z, m\in \mathbb Z_{\geq 0}$) by a formula
\[
\nu(B)=f(1_{B\cap A}).
\]
Here $1_X$ is the indicator of a set $X$ and since $B\cap A\subseteq \mathbb N$ the function $1_{B\cap A}$ is an element of $\ell^\infty$. Then $\nu$ can be extended uniquely to a Borel measure on $\mathbb Z_P$ and $\nu(X)\leq \mu(X)$ for any Borel set $X\subseteq \mathbb Z_P$ ($\mu$ is the Haar probability measure). Moreover, if the lower density of $A$ is positive, then $\nu$ is non-zero.
\end{lemma}

\begin{proof}
The last statement is immediate:
\[
\nu(\mathbb Z_P)=f(1_A)=f(C1_A)\geq \liminf_{n\to +\infty}C1_A(n)=\liminf_{n\to +\infty}\frac{\#A\cap[1,n]}{n}>0.
\]
For $C=a+P^m\mathbb Z_P$ we get
\[
\nu(C)\leq f(1_{a+P^m\mathbb Z})=\frac{1}{P^m}=\mu(a+P^m\mathbb Z_P)
\]
by strong invariance of $f$. Consequently, $\nu\leq \mu$ on $\mathcal C$. Therefore, it is also $\sigma$-additive on $\mathcal C$: if $C=\coprod_i C_i$, then
\[
\sum_{i=1}^{n}\nu(C_i)\leq \nu(C)\leq \sum_{i=1}^n \nu(C_i)+\sum_{i=n+1}^{+\infty}\mu(C_i)
\]
and taking limit $n\to +\infty$ we get $\sigma$-additivity, since the rightmost sum tends to $0$. Now, the Borel $\sigma$-algebra of $\mathbb Z_P$ is generated by clopen sets and by Carath\'eodory's extension theorem we get the desired result.
\end{proof}

Now take a semigroup $W$ generated by a balanced set $A_1,\ldots,A_r$ with $P=\prod\limits_{p\mid a_1\ldots a_r}p$. Let $B$ be a zero density subset of $\mathbb N$.

Now we are ready to prove Theorem 1.

\begin{proof}
Assume that $B\subseteq \mathbb N$ has density $0$, but the lower density of $W(B)$ is non-zero. Consider the measure $\nu=\nu_{W(B),f}$ from Lemma 2. The measure $\nu$ is non-zero and majorized by the Haar probability measure on $\mathbb Z_P$. What is left to prove is that $\nu$ is $W$-stationary. To do so, it is enough to prove that
\[
\nu(C)=\sum_{i=1}^r\frac{1}{a_i}(A_i)_*\nu(C)
\]
for basic sets $C=a+P^m\mathbb Z_P$, since they generate the Borel $\sigma$-algebra. Now,
\[
W(B)=\bigcup_{i=1}^r A_iW(B)\cup B,
\]
hence
\[
C\cap W(B)=\bigcup_{i=1}^r C\cap A_iW(B)\cup (C\cap B).
\]
Therefore,
\[
\nu(C)\leq \sum_{i=1}^r f(1_{C\cap A_iW(B)})+f(C\cap B).
\]
The last summand is zero, since $B$ has density $0$ and
\[
0\leq f(B)\leq \limsup_{n\to +\infty}\frac{\#B\cap [1,n]}{n}=0.
\]
On the other hand, $C\cap A_iW(B)=A_i(A_i^{-1}C\cap W(B))$ and due to strong invariance for any $X$ we have $f(1_{A_iX})=\frac{1}{a_i}f(1_X)$. Consequently,
\[
\nu\leq\sum_{i=1}^r a_i^-1(A_i)_*\nu.
\]
Since the two measures have the same total mass by the balance condition, the inequality must be an equality. Thus $\nu$ is $W$-stationary, completing the proof.
\end{proof}

Let us now prove Theorem 4. Its proof is fairly independent of previous considerations.
\begin{proof}[Proof of Theorem 4:]
Let $\omega(x)$ be the counting function of $W(B)$ and $b(x)$ be the counting function of $B$. By our assumptions, $b(x)=O(x^{1-c})$. Note that
\[
W(B)=\bigcup_{i=1}^r A_iW(B)\cup B.
\]
Therefore,
\[
\omega(x)\leq \sum_{i=1}^r \omega\left(\frac{x-b_i}{a_i}\right)+b(x)\leq \sum_{i=1}^r \omega\left(\frac{x}{a_i}\right)+b(x)
\]
by monotonicity. We introduce an auxiliary function $\beta(x)$ to absorb the error term $b(x)$. Define
\[
\beta(x)=x^{1-c/2}\sup_{y\geq x}\frac{b(y)}{y^{1-c/2}}.
\]
By construction, we have $\beta(x)\geq b(x)$, $\beta(x)=o(x)$ for $x\to +\infty$ and additionally for real $a>1$ we also have
\[
a^{1-c/2}\beta(x/a)=x^{1-c/2}\sup_{y\geq x/a}\frac{b(y)}{y^{1-c/2}}\geq x^{1-c/2}\sup_{y\geq x}\frac{b(y)}{y^{1-c/2}}=\beta(x), 
\]
which is useful for the next step. We now take a large enough constant $C$ and define $\psi(x)=\omega(x)+C\beta(x)$. Clearly,
\[
\underline{d}=\liminf_{x\to +\infty}\frac{\psi(x)}{x} \text{ and }\overline{d}=\limsup_{x\to +\infty}\frac{\psi(x)}{x}.
\]
Plugging the formula $\omega(x)=\psi(x)-C \beta(x)$ into the inequality of $\omega(x)$ and estimating $b(x)\leq \beta(x)$, we obtain
\[
\psi(x)\leq \sum_{i=1}^r \psi\left(\frac{x}{a_i}\right)-C\sum_{i=1}^r\beta\left(\frac{x}{a_i}\right)+(C+1)\beta(x).
\]
Since $\beta(x/a)\geq a^{c/2-1}\beta(x)$ for $a>1$, in the second summand we get
\[
\sum_{i=1}^{r} \beta\left(\frac{x}{a_i}\right)\geq \sum_{i=1}^{r} a_i^{c/2-1} \beta(x)=L\beta(x)
\]
with $L>1$, since the set $A_1,\ldots,A_r$ is balanced. Taking $C=(L-1)^{-1}$ we conclude that
\[
 \psi(x)\leq \sum_{i=1}^r \psi\left(\frac{x}{a_i}\right).
\]
Fix $\varepsilon>0$, then by definition of $\underline{d}$, there is a number $x_0\geq 1$ with $\psi(x_0)\leq (\underline{d}+\varepsilon)x_0$. Let us prove that $\psi(x)\leq (\underline{d}+\varepsilon)\max a_i \cdot x$ for $x\geq x_0(\max a_i)^{-1}$. We prove by induction on $n$ that the desired inequality holds for $x\in[x_0M^{n-1},x_0M^n]$, where $M=\max a_i$. For $n=0$ the inequality is automatic: if $x_0M^{-1}\leq x\leq x_0$
\[
\psi(x)\leq \psi(x_0)\leq (\underline{d}+\varepsilon)x_0\leq (\underline{d}+\varepsilon)Mx.
\]
Suppose that we've proved our inequality for $n$ and it is false for some $x\in [x_0M^n,x_0M^{n+1}]$. Then it must also be false for some $x/a_i$ with $i\leq r$: otherwise
\[
\psi(x)\leq \sum_{i=1}^r (\underline{d}+\varepsilon)M\frac{x}{a_i}=(\underline{d}+\varepsilon)Mx.
\]
Since $r>1$ and all $a_i>1$, applying this procedure several times we will obtain the counterexample lying in the previous interval, which is a contradiction. This shows that
\[
\overline{d}\leq (\underline{d}+\varepsilon)M
\]
for arbitrary $\varepsilon>0$, which concludes the proof.
\end{proof}

\section{Non-existence of non-zero KR-measures}

In this section we prove the main result of the paper, namely Theorem 2. In contrast to the general analytic arguments of the preceding section, the proof of Theorem 2 is largely combinatorial. From now on, we assume that $W$ is generated by a set $A_1,\ldots,A_r$ satisfying the conditions of Theorem 2. Let us first show that $W$ is a free semigroup on $A_1,\ldots,A_r$, i.e. there are no non-trivial relations between the generators. This result is the same as in Kolpakov-Talambutsa \cite{KolpakovTalambutsa2022}, but we give the proof for the sake of completeness.

\begin{lemma}
If a balanced semigroup $W=\langle A_1,\ldots, A_r\rangle$ satisfies the condition 2 of Theorem 2, then it is free.
\end{lemma}
\begin{proof}
Suppose, for contradiction, that $W$ is not free. There exist distinct words $w_1,w_2$, representing the same affine map, whose first letters differ and whose last letters differ. Then we can assume they are of the same length, since $w_1,w_2$ can be replaced by $w_1w_2, w_2w_1$. Next, from condition 2 we see that there is a compact set $K\subset \mathbb R$ such that $K$ has non-empty interior, hence also non-zero Lebesgue measure $\lambda(K)$, and for any $x\in K$ there is $i\leq r$ with $A_ix\in K$, so we have an inclusion
\[
K\subset \bigcup_{i=1}^r A_i^{-1}(K).
\]
Let $\ell$ be the shared length of our words $w_1,w_2$, we iterate the inclusion $\ell$ times to get
\[
K\subset \bigcup_{\#w=\ell} w^{-1}(K).
\]
Here $\#w$ denotes the length of $w$ in the alphabet $A_1,\ldots,A_r$. Since the set of generators is balanced, we have
\[
\sum_{\#w=\ell} \frac{1}{a_w}=\left(\sum_{i=1}^r \frac{1}{a_i}\right)^\ell=1.
\]
Notice that $\lambda(w^{-1}K)=\frac{1}{a_w}\lambda(K)$ for all $w$. Since $w_1$ and $w_2$ coincide as maps, but not as words, we get
\[
\lambda(K)\leq \lambda\left(\bigcup_{\#w=\ell} w^{-1}K\right)=\lambda\left(\bigcup_{\substack{\#w=\ell\\ w\neq w_1}} w^{-1}K\right)\leq \left(1-\frac{1}{a_{w_1}}\right)\lambda(K)<\lambda(K),
\]
which is a contradiction.
\end{proof}

Assume now that $\nu$ is a KR-measure for $W$. We also normalize it to be a probability measure.
\begin{lemma}
Suppose that $W'\subset W$ is primitive: for $w_1, w_2\in W'$ with $w_1\neq w_2$ the word $w_1$ is not a prefix for $w_2$. Then
\[
\nu\geq \sum_{w\in W'}\frac{1}{a_w}w_*\nu.
\]
\end{lemma}
\begin{proof}
We apply the $W$-stationarity relation
\[
\nu=\sum_{i=1}^r \frac{1}{a_{A_i}}(A_i)_*\nu
\]
and iterate it for each summand. Every new summand arising at any step has a form $\frac{1}{a_w}w_*\nu$. If we reach $w\in W'$, we stop the iteration, otherwise we continue. Since $W'$ is primitive, each $w\in W'$ eventually appears in our formula: indeed, if we assume the contrary, we must ``block'' a word $w_1$ with some $w_2$, i.e. we stop at $w_2$ and cannot reach $w_1$ anymore. This happens precisely when $w_2$ is a prefix for $w_1$.
\end{proof}
Now we fix some notations. For $w\in W$, $wx=a_wx+b_w$, $n(w)=v_p(a_w)$ --- the $p$-adic valuation of the leading coefficient and $m(w)=v_q(a_w)$. From the arithmetic condition for $W$ we see that
\[
a_w=p^{n(w)}q^{m(w)}.
\]
For a fixed $n\geq 0$ denote by $W_n$ the set of all $w\in W$ such that $n(w)=n$ and the last letter $A_j$ of $w$ has $n(A_j)>0$, i.e. $a_{A_j}$ is divisible by $p$. Clearly, $W_n$ is a primitive set. Indeed, if $w_1,w_2\in W_n$ and $w_2$ is a prefix for $w_1$, then $w_1=w_2wA_j$ for some $w$, but then $n(w_1)\geq n(w_2)+n(A_j)>n(w_2)$, which is a contradiction.

\begin{definition}
We say that the element $w\in W$ is nice if $b_w\equiv 0 \pmod{p^{n(w)}}$. We call a nice $w$ reduced if it does not begin with $A_1$.
\end{definition}
In other words, $w\in W$ is nice if it sends the whole $\mathbb Z_P=\mathbb Z_{pq}$ in the smallest possible $p$-adic neighbourhood of $0$. And a nice $w$ is reduced if it is nice for ``non-trivial reasons'': clearly for any nice $w$ the element $A_1w$ is also nice.
\begin{lemma}
For any $n\geq 0$ we have
\[
\nu(p^n\mathbb Z_P)\geq \sum_{\substack{w\in W_n\\ w\text{ is nice}}}p^{-n}q^{-m(w)}.
\]
\end{lemma}
\begin{proof}
Applying Lemma 4, we see that
\[
\nu(p^n\mathbb Z_P)\geq \sum_{w\in W_n}\frac{1}{a_w}w_*\nu(p^n\mathbb Z_P).
\]
Now, $w_*\nu(p^n\mathbb Z_P)=\nu(w^{-1}p^n\mathbb Z_P)$. But $w(\mathbb Z_P)\subset b_w+p^n\mathbb Z_P$, so the last measure is $1$ if $w$ is nice and $0$ otherwise, which concludes the proof.
\end{proof}

To obtain a contradiction, it remains to construct sufficiently many nice words. Our real and discrete dynamics conditions help with that. Let $K$ be the compact set from condition 2 and $I\subset K$ be an interval contained in $K$. For a fixed $m$, consider the set $Q_m$ of all rationals of the form $h/q^m$ with $h\in \mathbb Z$ and $h/q^m\in K$. Clearly, for large enough $m$ we have $|Q_m|\gg q^m$, because there are many rationals in an interval. Let us also define $R_m$ to be the subset of $Q_m$: the numbers $h/q^m$ with $h\equiv 0 \pmod {p^N}$. Similarly, $|R_m|\gg q^m$.

\begin{lemma}
For every non-integer $s\in p^{-N}R_m$ there is a reduced nice element $u=u(s)$ such that $us=0$ and the last letter of $u$ is $A_1$. The map $u:p^{-N}R_m\setminus \mathbb Z\to W$ is injective.
\end{lemma}
\begin{proof}
 Let $F:K\to K$ be a map such that for any $x\in K$ there is $i\leq r$ with $F(x)=A_ix$. We claim that for any non-integral $s\in Q_m$ there is an integer $S$ such that $F^{\circ S}(s)=a/q^j$ for some $(a,q)=1$ and $1\leq j\leq v$. 

We argue by induction over $m$. The cases $1\leq m\leq v$ are trivial. What is left to show is that for $m>v$ after several steps $F^{\circ S_1}(s)$ must lie in $Q_{m-1}$ and still be non-integral. Since all $A_i$ with $i>1$ have leading coefficients divisible by $q$ but not by $q^{v+1}$, the only case this could be false is if $F^{\circ n}(s)=A_1F^{\circ (n-1)}(s)=p^NF^{\circ(n-1)}(s)=\ldots=p^{Nn}s$ for all $n$. But this is impossible: $s\neq 0$, so $p^{Nn}|s|\to +\infty$ for $n\to +\infty$, which contradicts the compactness of $K$. Hence, the claim follows: there is an $S$ such that $F^{\circ S}(s)\in Q_v\setminus\mathbb Z$. 

Since every iteration of $F$ is some $A_i$, we have $A_{i_1}\ldots A_{i_S}s=a/q^j$ for some $1\leq j\leq v$ and $(a,q)=1$. Due to discrete dynamics condition, there is an element $w\in W$ such that $wA_{i_1}\ldots A_{i_S}s=0$.  Set $u'=wA_{i_1}\ldots A_{i_S}$. If $u'$ is not reduced, then $u'=A_1^lu$ for some integer $l\geq 1$ and reduced $u$. Since $0$ is the fixed point of $A_1$, we get $us=0$. Take now $s\in p^{-N}R_m$. By definition, $A_1s \in R_m$ and the denominator of $s$ is a power of $q$. Take then $u=u'A_1$, where $u'$ is the element corresponding to $p^Ns$.

Let us show that $u$ is nice. By the above observations, $u\in W$ is nice if and only if $u(\mathbb Z_p)\subset p^{n(u)}\mathbb Z_p$, which holds if and only if $u(\mathbb Z_p)\cap p^{n(u)}\mathbb Z_p\neq \emptyset$. This intersection is, indeed, non-empty, because by construction it contains $0$, since $s\in \mathbb Z_p$. Injectivity follows from the fact that the equation $ux=0$ has at most one solution.
\end{proof}

Let us now finish the proof of Theorem 2.

\begin{proof}
For large enough $m$, we have $|R_m\setminus R_{m-1}|\gg q^m$ (there are many rational points with a given denominator in an interval). By construction, if $p^Ns\in R_m\setminus R_{m-1}$, then $m(u(s))=m+O(1)$, because the construction of $u(s)$ involves applying our generators in some order until the denominator cancels out enough (which adds $m+O(1)$ to $v_q(r)=-m)$) and then reaching $0$ in $O(1)$ steps. Therefore,
\[
\sum_{s\in p^{-N}(R_m\setminus R_{m-1})}q^{-m(u(s))}\gg 1
\]
for $m\geq m_0$. Taking $m=m_0+M$ for a large constant $M$, we get finitely many nice reduced elements $u(s)$ for $s\in p^{-N}R_m$. Let $L=L(M)$ be the maximal value of $n(u)$ for elements of this form. Take $n>L(M)$, then for each $s$ we can find $t\in \mathbb Z_{\geq 0}$ such that $A_1^tu(s)\in W_{n}\cup W_{n+1}\cup \ldots \cup W_{n+N-1}$. All these words are nice, distinct from each other due to injectivity of $u$. For each denominator level $j$, the corresponding words have weights $q^{-m(w)}\asymp q^{-j}$, and their total contribution is bounded below by a positive constant independent of $j$. So we get
\[
\sum_{i=0}^{N-1}\sum_{\substack{w\in W_{n+i}\\ w\text{ is nice}}}q^{-m(w)}\gg M.
\]
Applying Lemma 5 we see that
\[
\sum_{i=0}^{N-1}p^{n+i}\nu(p^{n+i}\mathbb Z_P)\gg M.
\]
On the other hand, $\nu$ is bounded by $C\mu$ for some constant $C$, so the left-hand side of this inequality is $O(N)$. This gives $M\ll N$, with an implied constant independent of $M$. Since $N$ is fixed and $M$ is arbitrary, we obtain a contradiction.
\end{proof}
\begin{remark}
A slightly more careful proof of Theorem 2 would result in a quantitative bound
\[
\nu(p^n\mathbb Z_P)\gg \sqrt{n}p^{-n},
\]
but this is not necessary for our considerations.
\end{remark}

\section{Klarner's six}

Here we prove our Theorem 3. Note first that the existence of KR-measure is not changed after conjugation. To be more precise, suppose that $Ux=ax+b$ is an invertible map on $\mathbb Z_P$, i.e. $b\in \mathbb Z_P$ and $a\in \mathbb Z_P^*$. If semigroups $W,W'$ have the same $P$ and their generators are conjugated via $U$, then the existence of KR-measure $\nu$ for $W$ implies the same for $W'$ and vice versa. Indeed, suppose that $A_1,\ldots,A_r$ are generators of $W$ and $A_1',\ldots,A_r'$ are generators of $W'$ and $A_i'=UA_iU^{-1}$. Then if $\nu$ is a KR-measure for $W$, the measure
\[
\nu'=U_*\nu
\]
is a KR-measure for $W'$. First, if $\nu=f\mu$, where $f$ is the bounded density, then $\nu'=(f\circ U^{-1})\mu$ and $f(U^{-1}x)$ is also clearly bounded. Next,
\[
\nu=\sum_i a_{A_i}^{-1} (A_i)_*\nu.
\]
Applying $U_*$ to both sides we see
\[
\nu'=U_*\nu=\sum_i a_{A_i}^{-1} U_*(A_i)_*\nu=\sum_{i} a_{A'_i}^{-1} (A_i')_*\nu',
\]
as needed.

Six Klarner's sequences are divided into three conjugacy classes: $\langle 2x,\,3x+2,\,6x+3\rangle$ is conjugate to $\langle 2x+1,\,3x,\,6x+2\rangle$, here $Ux=-x-1$; $\langle 2x,\,3x+3,\,6x+10\rangle$ is conjugate to  $\langle 2x+2,\,3x+1,\,6x\rangle$ via $Ux=-x-2$; finally, 
$\langle 2x+2,\,3x,\,6x+15\rangle$ is conjugate to  $\langle 2x+1,\,3x+6,\,6x\rangle$, here $Ux=-x-3$. Note also that these conjugacy relations would not make sense in the positive integer setting, since none of them map $\mathbb N$ to itself. Due to our observations, it is enough to work with $\langle 2x,\,3x+2,\,6x+3\rangle$, $\langle 2x,\,3x+3,\,6x+10\rangle$ and $\langle 2x+2,\,3x,\,6x+15\rangle$. Let us denote these $W_1,W_2,W_3$ and corresponding generators by $R_i,S_i$ and $T_i$. Hence, for example, $S_2=3x+3$ and $R_3=2x+2$.

We now apply Theorem 2 to $W_i$. For $W_1,W_2$ we take $p=2, q=3, N=1$ and for $W_3$ we swap $2$ and $3$: $p=3, q=2, N=1$. Next, we define $K_1,K_2,K_3$ and corresponding maps. Here $K_1=[-1,0]$, $K_2=[-2,0]$ and $K_3=[-3,0]$. The maps $F_i:K_i\to K_i$:

\[
F_1(x)=\begin{cases}
R_1x=2x \text{ for }x\in[-1/2,0]\\
T_1x=6x+3 \text{ for }x\in [-2/3,-1/2)\\
S_1x=3x+2\text{ for }x\in [-1,-2/3)
\end{cases}
\]
\[
F_2(x)=\begin{cases}
R_2x=2x \text{ for }x\in[-1,0]\\
S_2x=3x+3 \text{ for }x\in [-5/3,-1)\\
T_2x=6x+10 \text{ for }x\in[-2,-5/3)
\end{cases}
\]
\[
F_3(x)=\begin{cases}
S_3x=3x \text{ for }x\in [-1,0]\\
R_3x=2x+2 \text{ for }x\in [-5/2,-1)\\
T_3x=6x+15 \text{ for }x\in [-3,-5/2)\\
\end{cases}
\]

These formulas mean that the arithmetic condition and the real dynamics condition hold for $W_i$. We verify the discrete dynamics condition by exhibiting the required rational orbits. For $i=1$ we need to look at $-2/3$ and $-1/3$. In this case we have
\[
-\frac13\underset{2x}\to -\frac23\underset{3x+2}\to 0.
\]
Further, for $W_2$ we consider $-5/3,-4/3,-2/3$ and $-1/3$. We obtain
\[
-\frac13\underset{2x}\to-\frac23\underset{2x}\to-\frac43\underset{3x+3}\to-1\underset{3x+3}\to0\text{ and }-\frac53\underset{6x+10}\to 0.
\]
Finally, for $W_3$ the denominators are equal to $2$ and we observe
\[
-\frac52 \underset{6x+15}\to0 \text{ and }-\frac12\underset{3x}\to -\frac32\underset{2x+2}\to-1\underset{2x+2}\to 0.
\]
Hence all the conditions of Theorem 2 are satisfied, which concludes the proof of Theorem 3.

\section{Conclusion}

We have proved that all six of Klarner's sequences have asymptotic density zero. More generally, our argument provides a criterion for the non-existence of stationary measures with bounded density, applicable to a wider class of affine semigroups. For example, one can take any two distinct primes $p,q$ and consider the maps $px, pqx+q,\ldots,pqx+pq-1$. They form a balanced set and $K=[-1,0]$ can be split into (almost) disjoint union of preimages of itself under these maps. Also, the rationals $-a/q$ with $0<a<q$ are mapped to $0$ by $x\mapsto p^{k}x\mapsto pq(p^kx)+p^{k+1}a$, where $k$ is the minimal one with $p^{k+1}a\geq q$, hence Theorem 2 is also applicable in this case.

 Previous work \cite{KonyaginProtasovTalambutsa2025} fully resolved the case when the maps have form $nx+a_1,\ldots,nx+a_n$. In this setting, Fourier analysis on $\mathbb Z_n$ allows one to compute the Fourier transform of the corresponding stationary measure explicitly. We are not aware yet of a Fourier-analytic proof of density zero for Klarner's sequences. It is likely that in this case some non-commutative input is needed.

The author thanks S. Konyagin, P. Kosenko, V. Protasov, A. Talambutsa and K. Shamazov for many fruitful discussions.

\end{document}